# Эволюционный вывод простейшей модели бимодального расщепления спроса на городские передвижения

Михаил Яковлевич Блинкин *mblinkin@hse.ru*
*(директор Института экономики транспорта и транспортной политики, НИУ ВШЭ),*
Александр Владимирович Гасников *gasnikov@yandex.ru*
*(доцент кафедры Математических основ управления ФУПМ МФТИ)*
Сергей Сергеевич Омельченко *sergey.omelchenko@phystech.edu*
*(магистр МФТИ)*
Инна Васильевна Усик *lavinija@mail.ru*
*(научный сотрудник БФУ им. И. Канта)*

**Аннотация**

В работе предлагается простейшая динамическая модель бимодального расщепления спроса на городские передвижения, удовлетворяемого, соответственно, поездкам на личных автомобилях и на общественном транспорте. Отличительной особенностью предлагаемого в статье подхода является учет эволюционной составляющей в механизме расщепления.

**Ключевые слова:** расщепление спроса на передвижения, закон Ципфа–Парето, принцип сжимающих отображений.

## 1. Введение

Расщепление спроса на ежедневные трудовые передвижения, удовлетворяемого, соответственно, поездкам на личных автомобилях и на общественном транспорте (*Modal Split*) базируется на стандартной (и вполне правдоподобной) гипотезе, согласно которой горожанин выбирает способ передвижения по критерию минимума обобщенной цены поездки (*Generalized Cost*), складывающейся из двух компонент:

 – денежных затрат на совершение поездки (*Monetary Cost*),

 – «немонетарных затрат» (*Non-Monetary Cost*), исчисляемых как произведение времени поездки на цену единицы времени горожанина (*Value of Time*), вообще говоря, сугубо индивидуальную для каждого.

Денежные затраты на совершение поездки трактуются в практике транспортного планирования как «траты из кармана» (*Out of Pocket Price*). Для общественного транспорта это плата за проезд на тех или иных маршрутах и видах транспорта. Для автомобильной поездки – траты на моторное топливо и парковочные платежи; условно-постоянные компоненты затрат (амортизация транспортного средства, налоги и т.п.) в учет обычно не берутся, так как они слабо влияют на ежедневные решения горожанина по выбору способа совершения поездки.

Способы оценки времени поездки на одно и то же расстояние различается для общественного транспорта и личного автомобиля принципиальным образом. Для общественного транспорта это время можно считать фиксированным: применительно к внеуличным видами транспорта это очевидный факт, для наземного транспорта, работающего в общем потоке транспортных средств, неизбежные задержки в движении, обусловленные плотностью трафика, учтены в маршрутных расписаниях. Для автомобильных поездок ключевым обстоятельством является зависимость времени поездки от плотности трафика, то есть в конечном итоге, от количества претендентов на пользование улично-дорожной сетью.

В данной статье будет предложена эволюционная модель поведения жителей простейшего города (с достаточно большим числом жителей), состоящего всего из двух условных районов: селитебного (спального) и делового (рабочего). У каждого горожанина-автомобилиста для совершения ежедневной трудовой поездки есть возможность воспользоваться как личным автомобилем, так и общей для всех линией (маршрутом) общественного транспорта. В отличие от многих других работ по данной тематике (см., например, [1]) в данной статье мы сделаем акцент не только на описание равновесного расщепления горожан по выбору способа передвижения (на личном транспорте или общественном), но и на том как в реальном времени может происходить процесс "нащупывания" этого равновесия. Таким образом, исследуется вопрос устойчивости равновесия.

## 2. Простейшая эволюционная модель расщепления

Итак, рассматривается город, состоящий из двух районов: спального и рабочего. Каждый день жители города (все они живут в спальном районе, а работают, естественно, в рабочем) ездят на работу. Каждый из них имеет личный автомобиль. Кроме того, в городе имеется развитая сеть общественного транспорта. Таким образом, каждый житель имеет две альтернативные возможности для ежедневных трудовых поездок: личный автомобиль и общественный транспорте.

Ежедневные потери пользователей личного транспорта, оценивающих единицу (минуту) своего времени в $p \geq 1$ рублей, могут быть рассчитаны следующим образом

$$A_p(x) = a + pT(x),$$

где $a > 0$ – характеризует постоянные затраты (цена топлива и т.п.), $x \in [0,1]$ – доля жителей города, использующих личный автомобиль, $T(x)$ – функция, характеризующая то, как пользователи транспортной сети оценивают свои временные затраты. Обычно (см., например, [1]) $T(x)$ выбирают вида BPR-функций[1], т.е. $T(x) = T_0 + \gamma x^4$.

Ежедневные потери пользователей общественного транспорта, оценивающих единицу (минуту) своего времени в $p \geq 1$ рублей, могут быть рассчитаны следующим образом

$$B_p(x) \equiv b_1 + pb_2,$$

где $b_1 > 0$ – характеризует постоянные затраты (цена билета и т.п.), а $b_2 > 0$ можно понимать как время, потерянное в пути. В отличие от личного транспорта, для общественного транспорта считается, что $b_2$ не зависит от $x$ (метро, электропоезда, выделенные полосы и т.п.).

Будем считать, что жителей в городе много. Они расслоены по тому, во сколько рублей каждый из них оценивает единицу своего времени (потерянного в пути). Обычно это колеблется от 1 руб/мин до 10 руб/мин (это максимальное значение будем обозначать далее через $p_{\max}$). Введем зависимость $x(p)$ – доля жителей города, оценивающих одну минуту своего времени не меньше чем в $p$ рублей. Эта зависимость естественным образом восстанавливается [3] из закона распределения населения по доходу Ципфа–Парето (см. приложение). Обычно эту зависимость считают степенной $x(p) = p^{-\eta}$, где $\eta$ выбирают из диапазона 1 – 2.

Наложим теперь физически правдоподобные ограничения. Во-первых, будем считать, что денежная цена автомобильной поездки (постоянные затраты) дороже, чем в случае общественного транспорта

$$b_1 < a. \qquad (1)$$

Это предположение не так очевидно, как это представляется на первый взгляд. К примеру, в Москве до введения платной парковки в 2013 году "Out of Pocket Price" автомобильной поездки была заметно ниже цены пересадочной поездки «метро + трамвай». Понятно, что в таких условиях общественный транспорт предпочитало в основном «безлошадное» население. Ограничение (1) автоматически выполняется немедленно после введения платной парковки даже по самому щадящему тарифу.

Во-вторых, будем считать, что «на автомобиле всегда быстрее»:

$$T(1) < b_2. \qquad (2)$$

Для справедливости этого предположения принципиально наличие поправочного коэффициента, связывающего времена, потерянные на личном и общественном транспорте (важно считать эти времена не равноценными). Понятно, что в условиях затора, практически неизбежного при $x$ близких к единице, время поездки на метро заведомо будет меньше, чем на автомобиле $T(1) > b_2$. Вопрос, однако, заключается в том, что время, проведенное, соответственно, в вагоне общественного транспорта и в собственном автомобиле, трудно считать равноценным. Другими словами, в данной работе под $b_2$ понимается время, потерянное в пути на общественном транспорте, приведенное к масштабу времени, потерянному на личном транспорте. Такое приведение увеличивает "физическое" время в виду различной комфортности перемещений.

---
[1] Функции данного типа были введены на основе обширных эмпирических обследований американским Bureau of Public Roads (BPR) и много лет применяются в работах по транспортному планированию и теории транспортного потока [2].

Добавим для аккуратности, что для горожан с наиболее высокой ценой времени обобщенная цена автомобильной поездки всегда ниже, чем на общественном транспорте:

$$a + p_{\max} T(1) < b_1 + p_{\max} b_2. \tag{3}$$

Предположим, наконец, что для горожан с наиболее низкой ценой времени дело обстоит прямо противоположным образом: даже при самой низкой загрузке улично-дорожной сети обобщенная цена поездки на общественном транспорте ниже, чем у автомобильной поездки:

$$a + T(0) > b_1 + b_2. \tag{4}$$

Определим зависимость $p(x)$, как корень уравнения $A_p(x) = B_p$, т.е.

$$p(x) = \frac{a - b_1}{b_2 - T(x)}.$$

При условиях (1) – (4) выписанная формула корректно определяет монотонную гладкую зависимость $p(x) \in (1, p_{\max})$ при $x \in [0,1]$.

Представим себе такую динамику (повторяющуюся из дня в день). Каждый житель в $(k+1)$-й день смотрит на то, какая доля жителей $x^k$ использовала личный автомобиль в $k$-й день. Считаем, что такая информация (статистика) по вчерашнему дню общедоступна (например, благодаря каким-нибудь интернет сервисам, скажем, Яндекс.Пробки). Исходя из этой информации каждый житель, оценивающий минуту своего времени в $p$ рублей, оценивает (экстраполируя ситуацию вчерашнего дня на день сегодняшний, за неимением точной информации о $x^{k+1}$) свои затраты от двух возможных альтернатив: $A_p(x^k)$ – личный автомобиль и $B_p$ – общественный транспорт. Мы считаем всех жителей рациональными, поэтому из двух альтернатив, каждый житель выбирает ту, которая приносит ему наименьшие затраты. Таким образом, происходит формирование $x^{k+1}$.

Из описанного выше ясно, что жители города в $(k+1)$-й день, оценивающие единицу своего времени в $p(x^k) < p \le p_{\max}$ рублей, предпочтут в этот день личный автомобиль, а жители, оценивающие единицу своего времени в $1 \leq p < p(x^k)$ рублей предпочтут в этот день общественный транспорт. Таким образом, в $(k+1)$-й день доля $x^{k+1} = x(p(x^k))$ жителей города использует (выберет) личный автомобиль.

Для того чтобы сформулировать основной результат, сделаем одно упрощающее предположение, которое позволит представить этот результат в более наглядной форме. Будем считать, что в зависимости $x(p) = p^{-\eta}$ параметр $\eta = 1$. Тогда, если

$$4\gamma < a - b_1, \qquad (5)$$

то $x(p(\cdot))$ – сжимающее преобразование отрезка $[0,1]$ в себя. Легко понять, что условия (1) – (5) совместны.

**Теорема 1.** *При $\eta = 1$ в условиях (1) – (5) описанная выше динамика сходится (вне зависимости от $x^0$) со скоростью геометрической прогрессии к неподвижной точке: $x^k \xrightarrow[k \to \infty]{} x^*$, которая определяется как единственная точка пересечения графиков $x(p)$ и $p(x)$ на плоскости $(p, x)$.*

**Доказательство.** Сформулированное в теореме 1 утверждение сразу следует из принципа неподвижной точки для сжимающих операторов [4]. Однако намного полезнее представляется продемонстрировать доказательство рисунком 1, проясняющим, как происходит "нащупывание" равновесия.

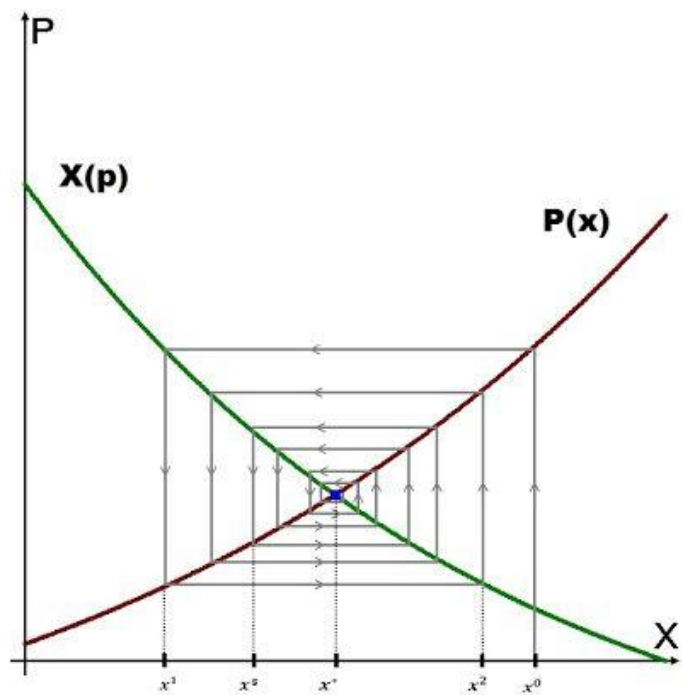

**Рис. 1**

Этот рисунок можно проинтерпретировать следующим образом. В начальный момент (в условный первый день) $x^0$ доля жителей города воспользовались личным транспортом. На следующий день всем известно $x^0$, исходя из этого числа, каждый оценивает каким образом ему сегодня добираться до работы (и обратно). В результате личным транспортом воспользуется $x^1 = x\left(p\left(x^0\right)\right)$ доля жителей города. Аналогично на следующий день личным транспортом воспользуется $x^2 = x\left(p\left(x^1\right)\right)$, на следующий $x^3 = x\left(p\left(x^2\right)\right)$ и т.д. Из рис. 1 видно, что такой итерационный процесс будет сходящимся к корню уравнения $x = x\left(p\left(x\right)\right)$, что можно понимать как точку пересечения графиков функции $x\left(p\right)$ и $p\left(x\right)$. Легко проверить, что если на рис. 1 выбрать $x^0$ левее $x^*$, то сходимость также будет иметь место. □

### 3. Численные эксперименты

Моделировался город численностью 1000 человек. Были выбраны следующие параметры $\gamma = 2$ минуты, $T_0 = 70$ минут, $b_2 = 75$ минут, $b_1 = 50$ рублей, $a = 60$ рублей. При таких параметрах процесс "нащупывания" равновесия сходился в среднем (случайно выбиралась точка $x^0 \in [0,1]$) за 3 – 4 дня (итерации). Были также проведены численные эксперименты и при других значениях параметров. При самых неблагоприятных значениях, используемых в численных экспериментах, сходимость была в среднем за 7 дней. Таким образом, получено модельное подтверждение известного из опыта факта, что выход (нащупывание) равновесия осуществляется крайне быстро (недели всегда хватает).

**Приложение: закон Ципфа–Парето и процесс Юла**

Далее мы постараемся пояснить возникновение закона Ципфа в тексте статьи. Для этого мы сделаем следующее упрощающее предположение: каждый человек оценивает единицу своего времени в сумму, которую он зарабатывает в единицу времени. Упрощая еще больше, предположим, что богатство каждого жителя прямо пропорционально тому, сколько он зарабатывает в единицу времени (считаем, что все работают одинаковое время – у всех нормированный рабочий день). Таким образом, нужно показать степенной характер распределения населения по доходу. Для этого рассмотрим следующую игрушечную модель.

В некотором городе живет неограниченно много жителей (изначально банкротов), которые могут участвовать в "освоении" монеток. База индукции: сначала выбирается один житель, он получает одну монетку. Шаг индукции: в $(k+1)$-й день выбирается очередной новый житель (отличный от $k$ уже выбранных), он получает возможность участвовать в разыгрывании монеток, монетка с вероятностью $\alpha < 1$ равновероятно отдается одному из этих $k+1$ жителей, а с вероятностью $1-\alpha$ эта монетка отдается одному из $k$ старых жителей с вероятностью пропорциональной тому, сколько у него уже есть монеток, т.е. по принципу "деньги к деньгам" (в моделях роста Интернета этот принцип называют "preferential attachment").

Обычно эту стохастическую динамику изучают в приближении среднего поля. В данном контексте это означает, что $n_s(t) \simeq E[n_s(t)]$, где $n_s(t)$ – количество жителей, у которых в ровно $s$ монеток на $t$-й день. Далее выписывают на $n_s(t)$ систему зацепляющихся обыкновенных дифференциальных уравнений. Автомодельное притягивающее решение этой системы ищут в виде $n_s(t) \sim x_s^* t$. После разрешения соответствующих уравнений получают степенной закон для зависимости $x_s^* \sim s^{-(1+1/(1-\alpha))}$. Такой подход применительно к моделям роста интернета (и изучения степенного закона для распределения степеней вершин) довольно часто сейчас встречается (в том числе и в учебной литературе). В частности, этот подход описан в обзоре [5], и рассматривается в книге [6]. По-сути, в этом приложении нами был описан процесс, возникающий в работах 20-х годов XX века по популяционной генетике, получивший названия процесса Юла (см., например, обзор [7]).

Описанная модель также восходит к работе конца XIX века В. Парето, в которой была предпринята попытка объяснить социальное неравенство и к работе Г. Ципфа конца 40-х годов XX века, в которой была отмечена важность степенных законов в "Природе". Эти законы для большей популярности иногда преподносят, как принцип Парето или принцип 80/20 (80 % результатов проистекают всего лишь из двадцати процентов причин) – такие пропорции отвечают $x_s^* \sim s^{-2.1}$ [7]. Приведем примеры (не совсем, правда, точные): 80 % научных результатов получили 20 % ученых, 80 % пива выпило 20 % людей и т.п. Сейчас много исследований во всем мире посвящено изучению возникновения в самых разных приложениях степенных законов (распределение городов по населению, коммерческих компаний по капитализации, автомобильных пробок по длине). Особенно бурный рост возник в связи с изучением роста больших сетей (экономических, социальных, интернета), см., например, книги и работы M.O. Jackson'a, в частности [8]. В России это направление также представлено, например, в Яндексе в отделе А.М. Райгородского [9] (Гречников–Остроумова–Рябченко–Самосват, Леонидов–Мусатов–Савватеев), в ИПМ РАН в группе А.В. Подлазова [10].